\newif\ifemsstyle
\IfFileExists{ems-journal.sty}{\emsstyletrue}{\emsstylefalse}
\ifemsstyle
  \documentclass{article}
  \usepackage[lang=american]{ems-journal}
\else
  \documentclass[11pt]{article}
  \usepackage[T1]{fontenc}
  \usepackage[utf8]{inputenc}
  \usepackage{lmodern}
  \usepackage[margin=1in]{geometry}
  \usepackage{amsmath,amssymb,amsthm}
  \usepackage{hyperref}
\fi
\usepackage{mathtools}
\usepackage{microtype}
\usepackage{xurl}
\numberwithin{equation}{section}
\theoremstyle{plain}
\newtheorem{theorem}{Theorem}[section]
\newtheorem{proposition}[theorem]{Proposition}
\newtheorem{lemma}[theorem]{Lemma}
\newtheorem{corollary}[theorem]{Corollary}
\theoremstyle{definition}
\newtheorem{remark}[theorem]{Remark}
\newcommand{\R}{\mathbb R}
\newcommand{\T}{\mathbb T}
\newcommand{\D}{\mathbb D}
\newcommand{\one}{\mathbf1}
\newcommand{\eps}{\varepsilon}
\DeclareMathOperator{\sgn}{sgn}
\DeclareMathOperator{\Sign}{Sign}
\DeclareMathOperator{\sech}{sech}
\DeclareMathOperator{\Dom}{Dom}
\DeclareMathOperator{\ran}{ran}
\allowdisplaybreaks[1]
\begin{document}

\hypersetup{hidelinks,pdftitle={Sharp weak-type (p,p) inequalities for Riesz transforms}}
\title{Sharp weak-type $(p,p)$ inequalities for Riesz transforms, $1\le p\le2$}
\ifemsstyle
\titlemark{Weak-type $(p,p)$ inequalities for Riesz transforms}

\emsauthor{1}{
  \givenname{Rodrigo}
  \surname{Ba\~nuelos}}{R.~Ba\~nuelos}
\Emsaffil{1}{
  \department{Department of Mathematics}
  \organisation{Purdue University}
  \zip{47907}
  \city{West Lafayette, IN}
  \country{USA}
  \affemail{banuelos@purdue.edu}}

\emsauthor{2}{
  \givenname{M\"{a}hdi}
  \surname{Hormozi}}{M.~Hormozi}
\Emsaffil{2}{
  \organisation{Beijing Institute of Mathematical Sciences and Applications}
  \zip{101408}
  \city{Beijing}
  \country{China}
  \affemail{hormozi@bimsa.cn}}

\emsauthor{3}{
  \givenname{Alexander}
  \surname{Volberg}}{A.~Volberg}
\Emsaffil{3}{
  \department{Department of Mathematics}
  \organisation{Michigan State University}
  \address{619 Red Cedar Road}
  \zip{48824}
  \city{East Lansing, MI}
  \country{USA}
  \affemail{volberg@math.msu.edu}}
\else
\author{Rodrigo Ba\~nuelos\thanks{Department of Mathematics,
Purdue University, West Lafayette, IN 47907, USA.
Email: \href{mailto:banuelos@purdue.edu}{banuelos@purdue.edu}.}
\and M\"{a}hdi Hormozi\thanks{Beijing Institute of Mathematical Sciences
and Applications, Beijing 101408, China.
Email: \href{mailto:hormozi@bimsa.cn}{hormozi@bimsa.cn}.}
\and Alexander Volberg\thanks{Department of Mathematics,
Michigan State University, 619 Red Cedar Road, East Lansing,
MI 48824, USA.
Email: \href{mailto:volberg@math.msu.edu}{volberg@math.msu.edu}.}}
\date{}
\maketitle
\fi

\begin{abstract}
For $1\le p\le2$, we obtain the sharp weak-type $(p,p)$ inequality
for each classical Riesz transform on real-valued functions on
$\mathbb R^n$, with the same optimal constant as for the Hilbert
transform. The proof combines the Bellman function
method of \cite{BHV} with Janakiraman's sharp strip
potential \cite{Janakiraman}.
\end{abstract}

\ifemsstyle
\classification[42A50, 60G44]{42B20}
\makeatletter
\let\ps@titlepage\ps@empty
\makeatother
\maketitle
\else
\fi

\section{Introduction}\label{sec:introduction}

Davis's theorem determines the best constant in the weak-type $(1,1)$
inequality for the Hilbert transform on the line and the conjugate-function
operator on the circle \cite{Davis}. Janakiraman
\cite[Theorems~1.1--1.2]{Janakiraman} identified the sharp constants for
$1\le p\le2$ for orthogonal harmonic functions and martingales under
differential subordination. His strip potential gives the corresponding
sharp inequalities for conjugate functions.

In \cite{BHV}, the authors used the Burkholder--Bellman function method
to prove that Davis's constant is also the sharp weak-type $(1,1)$
constant for each classical Riesz transform in every dimension.
Their argument combines weighted analytic positivity, regularization,
and variational confinement. This accompanying paper builds on that
approach, incorporating Janakiraman's potential to obtain the sharp
weak-type $(p,p)$ constants for $1<p\le2$.

The principal additional step is a version of the Davis functional
inequality with a prescribed global mean. It allows us to integrate
translated Davis functions with the compensation needed to recover
$|x|^p$. We prove this refinement and the compensated representation,
while quoting the analytic and variational inputs from \cite{BHV}.
The resulting integral inequality also retains a mean correction on
the torus. The constants and the underlying sharp strip potential are
those of \cite{Janakiraman}. The argument here establishes the corresponding
operator inequality for the classical Riesz transforms.

For $1\le p\le2$, we use the normalization of
\cite[Theorem~1.1 and (1.6)]{Janakiraman}:
\begin{equation}\label{eq:K-integral}
 K_p=
 \left[\frac1\pi\int_{-\infty}^{\infty}
 \frac{\left|\frac2\pi\log|t|\right|^p}{1+t^2}\,dt\right]^{-1}.
\end{equation}
Equivalently,
\begin{equation}\label{eq:K-beta}
 K_p=\frac{\pi^{p+1}}
 {2^{p+2}\Gamma(p+1)\beta(p+1)},
 \qquad
 \beta(s)=\sum_{k=0}^{\infty}\frac{(-1)^k}{(2k+1)^s}.
\end{equation}
We recall this evaluation below. In particular,
\[
 K_1=\frac{\pi^2}{8\beta(2)}=K_D,
 \qquad K_2=1.
\]
Here $K_p$ is the coefficient in the distribution inequality. With
the convention
\[
 \|v\|_{p,\infty}
 =\sup_{\lambda>0}\lambda\,\mu\{|v|>\lambda\}^{1/p},
\]
the corresponding operator norm is $K_p^{1/p}$.

On $\R^n$ we use Lebesgue measure and the Fourier convention
$\widehat f(\xi)=\int_{\R^n}e^{-ix\cdot\xi}f(x)\,dx$.
On $\T^n=(\R/(2\pi\mathbb Z))^n$ we use normalized Haar measure.
As in \cite{BHV},
\[
 R_j=-\partial_j(-\Delta)^{-1/2};
\]
its multiplier is $-i\xi_j/|\xi|$ on $\R^n$ and $-ik_j/|k|$ on
$\T^n$, with value zero at the origin.

\begin{theorem}\label{thm:main}
Let $1\le p\le2$, $n\ge1$, and $1\le j\le n$.
On either $\Omega=\R^n$ or $\Omega=\T^n$, every real-valued
$f\in L^p(\Omega)$ satisfies
\begin{equation}\label{eq:main-weak}
 \lambda^p\mu\{|R_jf|>\lambda\}
 \le K_p\|f\|_p^p,
 \qquad \lambda>0.
\end{equation}
The constant $K_p$ is optimal for every $n$ and $j$.
\end{theorem}

The case $p=1$ is \cite[Theorem~1.1]{BHV}, and the upper estimate at
$p=2$ follows from the $L^2$ contraction. Thus the new weak-type
extension concerns $1<p<2$. For $1<p\le2$, the proof also yields
the integral estimate in Proposition~\ref{prop:classical-B}, its
torus mean correction \eqref{eq:torus-mean}, and the refinement
depending only on $R_jf$ in Corollary~\ref{cor:refinement}.
These integral estimates have equality examples on the torus.

The finite-space proof uses a conservative transverse semigroup, and
the approximation in Section~\ref{sec:classical} preserves that
property. The conclusions concern the individual transforms $R_j$
on real-valued inputs.

\subsection{Evaluation of the constant}

Put $a_p=K_p^{-1}$. We include the calculation from
\cite[Section~4, p.~918]{Janakiraman} to fix the normalization.
Symmetry and the substitution $t=e^s$ in
\eqref{eq:K-integral} give
\begin{equation}\label{eq:a-sech}
 a_p=\left(\frac2\pi\right)^{p+1}
       \int_0^\infty s^p\sech s\,ds
     =\frac12\int_\R |x|^p\sech\left(\frac{\pi x}{2}\right)dx.
\end{equation}
Since $\sech s=2\sum_{k\ge0}(-1)^ke^{-(2k+1)s}$ for $s>0$,
\[
 \int_0^\infty s^p\sech s\,ds
 =2\Gamma(p+1)\sum_{k=0}^\infty
       \frac{(-1)^k}{(2k+1)^{p+1}}.
\]
Termwise integration is justified by
$\sum_{k\ge0}\int_0^\infty s^p e^{-(2k+1)s}\,ds<\infty$
for $p>0$. This proves \eqref{eq:K-beta}.

Section~\ref{sec:inputs} states the inputs from \cite{BHV}.
Section~\ref{sec:mean} proves the prescribed-mean refinement, and
Section~\ref{sec:superposition} derives the finite-space estimate by
compensated integration. Section~\ref{sec:classical} passes to the
classical transforms and states the integral refinements.
Section~\ref{sec:sharpness} proves sharpness in every fixed dimension
and records equality in the torus integral estimates.

\section{Analytic and variational inputs}\label{sec:inputs}

We use the notation and results below from \cite{BHV}. All references
to its equation and lemma numbers concern the version specified in
the bibliography.

\subsection{The Davis function}

Let $\Pi=\{x+iy:|y|<1\}$, where $x,y$ are the scalar
arguments of the Bellman functions. For $y\in(-1,1)$, let
$\mathsf P_y^\Pi(s)$ denote the density of the horizontal
coordinate of planar Brownian motion started at $(0,y)$,
at its first exit from $\Pi$. Here $y$ is the starting
height and $s\in\R$ is the horizontal exit coordinate.
Summing the Poisson kernels for the two boundary lines gives
\begin{equation}\label{eq:strip-density}
 \mathsf P_y^\Pi(s)=
 \frac{\cos(\pi y/2)\cosh(\pi s/2)}
 {2\bigl(\sinh^2(\pi s/2)+\cos^2(\pi y/2)\bigr)},
 \qquad |y|<1,\quad s\in\R.
\end{equation}
This is an even probability density, with
$\mathsf P_0^\Pi(s)=\tfrac12\sech(\pi s/2)$.
For Brownian motion started at $(x,y)$, the same density describes
the horizontal displacement before exit.

Let $u_1$ be the harmonic extension of $|x|$ to $\Pi$,
extended by $|x|$ for $|y|\ge1$, and put
\begin{equation}\label{eq:U1}
 u_1(x,y)=\int_\R |x+s|\mathsf P_y^\Pi(s)\,ds\quad(|y|<1),
 \qquad a_1=u_1(0,0),\qquad U_1=a_1-u_1.
\end{equation}
Thus $a_1=K_D^{-1}$ and $U_1$ is the function denoted by $U$ in
\cite{BHV}. By \cite[Section~2 and Lemma~3.3]{BHV}, it is continuous,
superharmonic, even in each variable, concave and one-Lipschitz in
its first variable, and
\begin{equation}\label{eq:U1-bounds}
 -|x|\le U_1(x,y)\le a_1-|x|,
 \qquad U_1(x,y)=a_1-|x|\quad(|y|\ge1).
\end{equation}
Moreover,
\[
 \omega(\delta):=\sup_{x,y,\,|h|\le\delta}
 |U_1(x,y+h)-U_1(x,y)|\longrightarrow0
 \quad(\delta\downarrow0).
\]
The analytic gradient in the strip satisfies
\begin{equation}\label{eq:F}
 F=u_{1,x}-iu_{1,y},\qquad F(0)=0,\qquad
 F'(z)=\sech(\pi z/2),\qquad \Re F'>0,
\end{equation}
and $|\Re F|\le1$; see \cite[Section~3.1]{BHV}.

\subsection{Finite conservative products}
Let $E$ be a finite space with weights $\nu_a>0$ summing to one,
and let $L\ge0$ be self-adjoint on $L^2(E,\nu)$.
Assume that $P_t=e^{-tL}$ is a conservative Markov semigroup.
Here $t\ge0$ is time, and $P_t$ acts on the state variable
$a\in E$; the parameters $y,s$ in \eqref{eq:strip-density}
are spatial coordinates. For $a,b\in E$, define the transition
probabilities by
\[
 p_t(a,b):=(P_t\mathbf1_{\{b\}})(a),
\]
so that, for every function $v:E\to\mathbb R$,
\[
 (P_tv)(a)=\sum_{b\in E}p_t(a,b)v(b).
\]
Positivity, conservativity, and self-adjointness give
\[
 p_t(a,b)\ge0,\qquad
 \sum_{b\in E}p_t(a,b)=1,\qquad
 \nu_a p_t(a,b)=\nu_b p_t(b,a).
\]
On
$X=\T\times E$, with $\mu=dm\,d\nu$, set
\begin{equation}\label{eq:operators}
 D=|D_x|,\qquad A=(-\partial_x^2+L)^{1/2},\qquad
 R=-\partial_xA^{-1},\qquad H=-\partial_xD^{-1}.
\end{equation}
Spectral inverses vanish on their kernels. Both $R$ and $H$ vanish
at horizontal frequency zero; $R^*=-R$ and $\|R\|_{2\to2}\le1$.
On functions with zero horizontal mean, put $T=AD^{-1}-I$.
The correction formula \cite[(3.10)]{BHV} is
\begin{equation}\label{eq:T}
 T=\frac1{2\sqrt\pi}\int_0^\infty H_t(I-P_t)t^{-3/2}\,dt,
 \qquad
 \widehat {H_tv}(k)=|k|^{-1}e^{-tk^2}\widehat v(k)\quad(k\ne0),
\end{equation}
where $H_t$ vanishes at zero and the integral converges in operator norm.
The operators $H_t=D^{-1}e^{-tD^2}$ act horizontally in each
state; $H$ in \eqref{eq:operators} denotes the Hilbert transform.
If $\bar f_a=\int_\T f_a\,dm$ denotes the horizontal mean
in state $a$,
\cite[(3.17)]{BHV} gives
\begin{equation}\label{eq:inverse}
 \|f-\bar f\|_2\le\sqrt{1+\|L\|}\,\|Rf\|_2,
 \qquad
 \|\bar f\|_2\le(\min_a\nu_a)^{-1/2}\|f\|_1.
\end{equation}
We use $W^{1/2,2}=W^{1/2,2}(\T;L^2(E,\nu))=\Dom(A^{1/2})$,
with equivalent norms. Globally Lipschitz maps act continuously by
composition on this space, by \cite[Lemma~3.7]{BHV}.

\subsection{Regularization and confinement}

Let $\zeta_\delta$ be an even smooth probability kernel supported in
$[-\delta,\delta]$, and let $U_{1,\delta}=U_1*_{y}\zeta_\delta$.
By \cite[Lemma~3.4]{BHV}, this preserves superharmonicity, evenness,
concavity and the one-Lipschitz bound in the first variable. Also
\begin{equation}\label{eq:smoothing}
 \|U_{1,\delta}-U_1\|_\infty\le\omega(\delta),\qquad
 U_{1,\delta}(x,y)=a_1-|x|\quad(|y|\ge1+\delta),
\end{equation}
and $\partial_yU_{1,\delta}$ is bounded and globally Lipschitz.
For $\eps>0$ put
\[
 J_{\delta,\eps}(f)=
 \int_XU_{1,\delta}(f,Rf)\,d\mu-\eps\|Rf\|_2^2,
 \qquad S_t=e^{-tD}.
\]
Here $S_t$ is the horizontal Poisson semigroup on $\T$,
acting separately in each state of $E$.

\begin{lemma}[{\cite[Lemma~3.5]{BHV}}]\label{lem:Poisson}
For real $f\in L^2(X)$ and $t>0$,
\begin{equation}\label{eq:Poisson}
 \begin{split}
 J_{\delta,\eps}(S_tf)\ge{}&J_{\delta,\eps}(f)
 +\eps\bigl(\|Rf\|_2^2-\|S_tRf\|_2^2\bigr)\\
 &-C_{\delta,L}t\bigl(\|f\|_2+\|f\|_2^2\bigr),
 \end{split}
\end{equation}
where the constant is independent of $\eps,t,f$.
\end{lemma}

The next criterion is extracted from the confinement calculation in
\cite[Section~3.3, (3.26)--(3.27)]{BHV}. That calculation uses only the
displayed hypotheses; it does not require $f$ to be an unconstrained
maximizer or $\rho$ to have zero mean.
Write $\Sign(r)=\{\sgn r\}$ for $r\ne0$ and $\Sign(0)=[-1,1]$.

\begin{lemma}[Contact criterion]\label{lem:confinement}
Let $\delta,\eps>0$, $f,g,\rho\in W^{1/2,2}$, $g=Rf$, and
$\eta=\partial_yU_{1,\delta}(f,g)$. Suppose
\begin{equation}\label{eq:contact-input}
 \begin{gathered}
 |\rho|\le1,\qquad
 \rho\in\Sign(f)\quad\hbox{on }\{|g|>1+\delta\},\\
 A\rho=\partial_x\eta-2\eps\partial_xg,
 \qquad Ag=-\partial_xf,
 \end{gathered}
\end{equation}
where the derivative identities hold in $W^{-1/2,2}$. Then
\begin{equation}\label{eq:quoted-confinement}
 |g|\le1+\delta\quad\hbox{almost everywhere}.
\end{equation}
\end{lemma}

\section{A Bellman estimate with prescribed mean}
\label{sec:mean}

We continue in the finite conservative probability setting
\eqref{eq:operators}. A bar denotes the horizontal mean in each
state, and $\mathfrak m=\int_X f\,d\mu$ denotes the global mean.
The superposition in Section~\ref{sec:superposition} will use the
translated inputs $f\pm r$. Retaining their means
$\mathfrak m\pm r$ supplies the compensation needed to integrate
the translated Davis inequalities.

\begin{lemma}\label{lem:center}
Let $f\in L^2(X;\R)$, put $g=Rf$ and
$\mathfrak m=\int_X f\,d\mu$, and assume $\|g\|_\infty<1$.
Define
\[
 h_0=-Hg,\qquad h=\mathfrak m+h_0,\qquad
 q=\Re F(h+ig).
\]
Then
\begin{equation}\label{eq:common-positivity}
 \langle q,Th_0\rangle\ge0,\qquad
 \bar q_a=F(\mathfrak m)\quad\hbox{for every }a\in E.
\end{equation}
\end{lemma}

\begin{proof}
In each state, $g_a$ has zero horizontal mean. Thus $h_a+ig_a$
is the boundary value of an analytic function $Z_a\in H^2(\D)$,
where $\D$ is the unit disk, with $Z_a(0)=\mathfrak m$.
Its imaginary part is the Poisson
extension of $g_a$, so $Z_a(\D)\subset\Pi$.
We adapt the weighted positivity argument of
\cite[Section~3.1]{BHV} to this common real center.

For two fibers put $Z=Z_a-Z_b$ and $Q=F(Z_a)-F(Z_b)$. Then
$Z(0)=Q(0)=0$ and
\[
 Q=BZ,\qquad
 B(z)=\int_0^1F'\bigl(Z_b(z)+\theta(Z_a(z)-Z_b(z))\bigr)\,d\theta,
 \qquad \Re B\ge0.
\]
The last assertion follows from convexity of the strip and
\eqref{eq:F}. We have $Z\in H^2$. Also $F(Z_a)\in H^2$:
its real part is bounded, and its value $F(\mathfrak m)$ at the
origin fixes the additive constant of its harmonic conjugate.
Consequently $Q\in H^2$.

Write $Z(z)=\sum_{j\ge1}Z_jz^j$ and $Q(z)=\sum_{j\ge1}Q_jz^j$.
For $n\ge1$, put $A_n=\tfrac12\Re\sum_{j=1}^nQ_j\overline{Z_j}$.
The coefficients through order $n$ of
$B(z)\sum_{j=1}^nZ_jz^j$ agree with those of $Q$. Hence, for $0<r<1$,
\[
 \frac12\int_\T\Re B(re^{ix})
 \left|\sum_{j=1}^nZ_jr^je^{ijx}\right|^2dm(x)
 =\frac12\Re\sum_{j=1}^nQ_j\overline{Z_j}r^{2j}\ge0.
\]
Let $r\uparrow1$ to obtain $A_n\ge0$.

Since $Z,Q\in H^2$, the coefficient products are absolutely summable.
For every bounded nonincreasing sequence $w_j\ge0$, Abel summation gives
\[
 \frac12\Re\sum_{j\ge1}w_jQ_j\overline{Z_j}
 =w_\infty A_\infty+
 \sum_{n\ge1}(w_n-w_{n+1})A_n\ge0.
\]
Apply this to $w_j=j^{-1}e^{-tj^2}$. Taking real parts of the boundary
functions gives
\[
 \langle q_a-q_b,H_t(h_{0,a}-h_{0,b})\rangle\ge0.
\]
Using the transition probabilities introduced above,
conservativity gives
\[
 \langle q,H_t(I-P_t)h_0\rangle_{L^2(X,\mu)}
 =
 \sum_{a,b\in E}\nu_a p_t(a,b)
 \langle q_a,H_t(h_{0,a}-h_{0,b})\rangle_{L^2(\mathbb T,m)}.
\]
Interchanging $a,b$, applying the detailed-balance identity
$\nu_a p_t(a,b)=\nu_b p_t(b,a)$, and averaging, we obtain
\[
 \langle q,H_t(I-P_t)h_0\rangle_{L^2(X,\mu)}
 =
 \frac12\sum_{a,b\in E}\nu_a p_t(a,b)
 \langle q_a-q_b,H_t(h_{0,a}-h_{0,b})\rangle_{L^2(\mathbb T,m)}
 \ge0.
\]
The inequality follows from $p_t(a,b)\ge0$ and the preceding
positivity estimate for each pair $a,b$.
Integrating over $t>0$ with the weight in \eqref{eq:T}
proves the first assertion of \eqref{eq:common-positivity}.
Finally, the mean-value property applied to $\Re(F\circ Z_a)$,
together with $Z_a(0)=\mathfrak m$, gives
\[
 \bar q_a=\Re F(Z_a(0))
         =\Re F(\mathfrak m)=F(\mathfrak m).
\]
\end{proof}

\begin{lemma}\label{lem:mean-strip}
If $f\in L^2(X;\R)$ and $|Rf|\le1$ almost everywhere, then
\begin{equation}\label{eq:mean-strip}
 \int_X U_1(f,Rf)\,d\mu\le U_1(\mathfrak m,0),
 \qquad \mathfrak m=\int_X f\,d\mu.
\end{equation}
\end{lemma}

\begin{proof}
First assume $\|g\|_\infty<1$, where $g=Rf$, and use the notation
of Lemma~\ref{lem:center}. Each $h_a+ig_a$ is the boundary value
of an analytic function with center $\mathfrak m$. Harmonicity of
$U_1$ inside the strip gives
\[
 \int_\T U_1(h_a,g_a)\,dm=U_1(\mathfrak m,0).
\]
This follows first on smaller disks and then at the boundary from
radial $L^2$ convergence, the Lipschitz bound in the first variable,
and uniform continuity in the second.

Since $\partial_xU_1(h,g)=-q$, concavity gives
\[
 U_1(f,g)\le U_1(h,g)-q(f-h).
\]
The multiplier identity used in \cite{BHV} is
$f-\bar f=(I+T)h_0$. Thus
\[
 f-h=Th_0+\bar f-\mathfrak m.
\]
By Lemma~\ref{lem:center}, the first term has nonnegative pairing
with $q$, while the common fiber mean of $q$ gives
\[
 \int_Xq(\bar f-\mathfrak m)\,d\mu
 =F(\mathfrak m)\sum_a\nu_a(\bar f_a-\mathfrak m)=0.
\]
Integration proves \eqref{eq:mean-strip} in the strict strip.
For the closed strip, apply this to $rf$, $0<r<1$.
Its mean is $r\mathfrak m$, so the right side is
$U_1(r\mathfrak m,0)$. Let $r\uparrow1$, using
\eqref{eq:U1-bounds} and continuity.
\end{proof}

The next proposition strengthens \cite[(3.1)]{BHV} by retaining
the global mean and removes the strip restriction.

\begin{proposition}\label{prop:mean}
For every real $f\in L^2(X)$,
\begin{equation}\label{eq:mean}
 \int_X U_1(f,Rf)\,d\mu
 \le U_1\left(\int_X f\,d\mu,0\right).
\end{equation}
\end{proposition}

\begin{proof}
Fix $\mathfrak m\in\R$ and maximize $J_{\delta,\eps}$ over the
closed affine subspace
\[
 \mathcal A_{\mathfrak m}
 =\left\{v\in L^2(X;\R):\int_Xv\,d\mu=\mathfrak m\right\}.
\]
We follow \cite[Sections~3.2--3.3]{BHV}, keeping this constraint
throughout the regularization and confinement argument.

\emph{Existence and regularity.}
The constant function $\mathfrak m$ is admissible,
$R\mathfrak m=0$, and
\[
 J_{\delta,\eps}(\mathfrak m)
 =U_{1,\delta}(\mathfrak m,0)>-\infty.
\]
This finite lower bound, the estimate
\[
 J_{\delta,\eps}(v)\le a_1-\|v\|_1-\eps\|Rv\|_2^2,
\]
and the inverse bounds \eqref{eq:inverse} show that maximizing
sequences are bounded in $L^2$. The affine subspace is weakly
closed, and horizontal Poisson smoothing preserves its constraint.
If $v_n\rightharpoonup v$ in $L^2$, compactness of $S_t$ gives
$S_tv_n\to S_tv$ strongly for each fixed $t>0$. Strong continuity
of $J_{\delta,\eps}$ and \eqref{eq:Poisson}, with the nonnegative
energy gain discarded, imply
\[
 \limsup_nJ_{\delta,\eps}(v_n)
 \le J_{\delta,\eps}(S_tv)+Ct,
\]
where $C$ is uniform for the bounded sequence. Letting $t\downarrow0$
proves weak upper semicontinuity. A maximizer $f_*$ therefore exists.

Write $g_*=Rf_*$. Since $S_tf_*\in\mathcal A_{\mathfrak m}$,
maximality and \eqref{eq:Poisson} imply
\[
 \eps\bigl(\|g_*\|_2^2-\|S_tg_*\|_2^2\bigr)\le C_{f_*}t.
\]
Fourier expansion and Fatou's lemma give
\[
 2\sum_{k\in\mathbb Z}|k|\|\widehat g_*(k)\|_{L^2(E)}^2
 \le\liminf_{t\downarrow0}
 \frac{\|g_*\|_2^2-\|S_tg_*\|_2^2}{t}
 \le C_{f_*}/\eps.
\]
Thus $g_*\in W^{1/2,2}$. The inverse multiplier bound, now with
the factor $|k|$ in the Fourier sum, gives the same regularity for
the nonzero horizontal frequencies of $f_*$. Its horizontal mean
lies in a finite-dimensional space, so $f_*\in W^{1/2,2}$.

\emph{The constrained Euler relation.}
Set $\eta=\partial_yU_{1,\delta}(f_*,g_*)$ and freeze the second
variable in the continuous convex functional
\[
 \mathcal C(v)=-\int_XU_{1,\delta}(v,g_*)\,d\mu.
\]
For $h\in L^2(X;\R)$ of global mean zero and $t\downarrow0$,
the Lipschitz bound for $\partial_yU_{1,\delta}$ gives
\[
 \begin{split}
 0&\ge J_{\delta,\eps}(f_*+th)-J_{\delta,\eps}(f_*)\\
  &=-\mathcal C(f_*+th)+\mathcal C(f_*)
    +t\langle\eta-2\eps g_*,Rh\rangle+o(t).
 \end{split}
\]
Thus the directional derivative of $\mathcal C$ is at least
$\langle w,h\rangle$, where
\[
 w=R^*(\eta-2\eps g_*)=-R(\eta-2\eps g_*).
\]
Convexity gives
\[
 \mathcal C(f_*+h)\ge\mathcal C(f_*)+\langle w,h\rangle
 \quad\hbox{when }\int_Xh\,d\mu=0.
\]
The normal space to $\mathcal A_{\mathfrak m}$ is the span of
$\one$. The subgradient sum rule, applicable because $\mathcal C$
is continuous, therefore gives a scalar $\kappa$ such that
\begin{equation}\label{eq:rho}
 \rho=w+\kappa\one\in\partial\mathcal C(f_*).
\end{equation}
Varying on measurable subsets identifies the subgradient pointwise:
\begin{equation}\label{eq:rho-contact}
 |\rho|\le1,\qquad
 \rho\in\Sign(f_*)\quad\hbox{on }\{|g_*|>1+\delta\}.
\end{equation}
Since $R$ has zero horizontal mean, $\kappa=\int_X\rho\,d\mu$
and hence $|\kappa|\le1$.

Lipschitz composition on $W^{1/2,2}$ gives
$\eta\in W^{1/2,2}$, and boundedness of $R$ on this space gives
$\rho\in W^{1/2,2}$. Conservativity implies $L\one=0$ and
$A\one=0$. Applying $A$ to \eqref{eq:rho} therefore gives
\[
 A\rho=\partial_x\eta-2\eps\partial_xg_*,
 \qquad Ag_*=-\partial_xf_*.
\]
Together with \eqref{eq:rho-contact}, these are the inputs
\eqref{eq:contact-input} to Lemma~\ref{lem:confinement}, so
\begin{equation}\label{eq:confined}
 |g_*|\le1+\delta\quad\hbox{almost everywhere}.
\end{equation}
The contact calculation uses the zero horizontal mean of $g_*=Rf_*$;
it imposes no zero-mean condition on $\rho$. Thus the scalar
multiplier in \eqref{eq:rho} is compatible with that calculation.

\emph{Removal of the regularization.}
Evenness and concavity in the first variable, uniform continuity
in the second, and \eqref{eq:confined} give
\[
 \begin{split}
 J_{\delta,\eps}(f_*)
 &\le\int_XU_1(f_*,g_*)\,d\mu+\omega(\delta)\\
 &\le\int_X U_1\left(\frac{f_*}{1+\delta},
                         \frac{g_*}{1+\delta}\right)d\mu
       +2\omega(\delta)\\
 &\le U_1\left(\frac{\mathfrak m}{1+\delta},0\right)
       +2\omega(\delta).
 \end{split}
\]
The last line is Lemma~\ref{lem:mean-strip}; the rescaled input has
mean $\mathfrak m/(1+\delta)$. For any fixed
$f\in\mathcal A_{\mathfrak m}$, maximality and \eqref{eq:smoothing}
imply
\[
 \int_XU_1(f,Rf)\,d\mu
 \le U_1\left(\frac{\mathfrak m}{1+\delta},0\right)
      +3\omega(\delta)+\eps\|Rf\|_2^2.
\]
Let $\eps\downarrow0$ and then $\delta\downarrow0$ to obtain
\eqref{eq:mean}.
\end{proof}

\section{Janakiraman's potential and compensated translates}
\label{sec:superposition}

For $1<p\le2$, let $u_p$ be the harmonic extension of $|x|^p$
to $\Pi$, with polynomial growth, and set $u_p(x,y)=|x|^p$
for $|y|\ge1$. Thus
\begin{equation}\label{eq:up}
 u_p(x,y)=\int_\R |x+s|^p\mathsf P_y^\Pi(s)\,ds
 \quad (|y|<1),\qquad u_p(0,0)=a_p.
\end{equation}
Define
\begin{equation}\label{eq:Up-Bp}
 U_p(x,y)=a_p-u_p(x,y),\qquad
 B_p(x,y)=U_p(x,y)+|x|^p.
\end{equation}
The function $u_p$ is Janakiraman's $W$, given in
\cite[pp.~912--913, (2.1)--(2.2)]{Janakiraman}, and
$a_p^{-1}U_p=1-K_pu_p$ is his superharmonic majorant $U$
\cite[p.~913]{Janakiraman}. We shall also use
\begin{equation}\label{eq:phi-definition}
 \Phi_p(y)=a_p-u_p(0,y).
\end{equation}
The following lemma includes the obstacle property of
\cite[Lemma~2.1(a), p.~912; proof on p.~916]{Janakiraman}.
We give a direct proof that also identifies a lower bound depending
only on the second variable.

\begin{lemma}\label{lem:obstacle}
The functions $u_p,U_p,B_p,\Phi_p$ are continuous. For $1<p\le2$,
\begin{equation}\label{eq:B-bounds}
 0\le\Phi_p(y)\le B_p(x,y)\le a_p,\qquad
 \Phi_p(y)=B_p(x,y)=a_p\quad (|y|\ge1).
\end{equation}
Moreover, $B_p(0,0)=\Phi_p(0)=0$ and $U_p(\mathfrak m,0)\le0$
for every $\mathfrak m\in\R$.
\end{lemma}

\begin{proof}
The strip densities in \eqref{eq:up} have exponential tails and
concentrate at zero as $y\to\pm1$. The representation therefore
gives continuity, including convergence to the boundary data.
Jensen's inequality gives $u_p(x,y)\ge|x|^p$ in the strip.

We use the scalar inequality
\begin{equation}\label{eq:symmetric-power}
 \frac{|x+s|^p+|x-s|^p}{2}\le |x|^p+|s|^p,
 \qquad 1\le p\le2.
\end{equation}
For $1<p\le2$, symmetry and homogeneity reduce this to
$(1+s)^p+(1-s)^p\le2+2s^p$ for $0\le s\le1$.
Writing $q=p-1\in(0,1]$, the derivative of the left side minus
the right side is
\[
 p\bigl[(1+s)^q-(1-s)^q-2s^q\bigr]
 \le p\bigl[(2s)^q-2s^q\bigr]\le0,
\]
by subadditivity of $r^q$; its value at zero is zero.
The case $p=1$ is immediate.

Symmetry of $\mathsf P_y^\Pi$ and \eqref{eq:symmetric-power} imply
\begin{equation}\label{eq:phi-comparison}
 u_p(x,y)-|x|^p\le u_p(0,y),\qquad
 B_p(x,y)\ge\Phi_p(y).
\end{equation}
The function $u_p$ is convex in $x$, harmonic inside the strip,
and even in $y$. Hence $u_{p,yy}=-u_{p,xx}\le0$, so
$0\le u_p(0,y)\le u_p(0,0)=a_p$.
This proves \eqref{eq:B-bounds} in the strip; outside it the
assertions follow from the definition. Finally,
$\mathfrak m\mapsto u_p(\mathfrak m,0)$ is even and convex,
so its minimum is $a_p$, proving $U_p(\mathfrak m,0)\le0$.
\end{proof}

Since $B_p(0,y)=\Phi_p(y)$, the function $\Phi_p$ is the largest
pointwise lower bound for $B_p(x,y)$ that depends only on $y$.

\begin{lemma}\label{lem:superposition}
Let $1<p\le2$ and $c_p=p(p-1)/2$. For every $(x,y)\in\R^2$,
\begin{equation}\label{eq:superposition}
 \begin{split}
 U_p(x,y)=c_p\int_0^\infty
 \bigl[U_1(x+r,y)+U_1(x-r,y)-2U_1(r,0)\bigr]
 r^{p-2}\,dr.
 \end{split}
\end{equation}
Here $r>0$ is a translation parameter. The displayed
compensated integral is absolutely convergent.
\end{lemma}

\begin{proof}
For $x\ge0$, the function $|x+r|+|x-r|-2r$ equals $2(x-r)$
on $0\le r\le x$ and zero on $r\ge x$. Thus
\begin{equation}\label{eq:power-mixture}
 |x|^p=c_p\int_0^\infty
       (|x+r|+|x-r|-2r)r^{p-2}\,dr.
\end{equation}
The integrand is nonnegative. Apply \eqref{eq:power-mixture}
to $x+s$ and integrate against $\mathsf P_y^\Pi(s)\,ds$. Tonelli's theorem
gives, in the strip,
\begin{equation}\label{eq:up-mixture}
 u_p(x,y)=c_p\int_0^\infty
 [u_1(x+r,y)+u_1(x-r,y)-2r]r^{p-2}\,dr.
\end{equation}
Outside the strip this is \eqref{eq:power-mixture}.
At $(0,0)$ it gives
\[
 a_p=2c_p\int_0^\infty[u_1(r,0)-r]r^{p-2}\,dr.
\]
Subtracting these finite nonnegative integrals proves
\eqref{eq:superposition}. The absolute value of their difference
is bounded by their sum, which also proves absolute convergence.
\end{proof}

\begin{proposition}\label{prop:p-finite}
In the finite conservative probability setting, let $1<p\le2$
and $f\in L^p(X;\R)$, and put $\mathfrak m=\int_Xf\,d\mu$.
Then
\begin{equation}\label{eq:p-mean}
 \int_X U_p(f,Rf)\,d\mu\le U_p(\mathfrak m,0)\le0.
\end{equation}
Equivalently, for every $\lambda>0$,
\begin{equation}\label{eq:B-finite}
 \lambda^p\int_X B_p(f/\lambda,Rf/\lambda)\,d\mu
 \le\|f\|_p^p+\lambda^pU_p(\mathfrak m/\lambda,0)
 \le\|f\|_p^p.
\end{equation}
For inputs outside $L^2$, $R$ denotes its extension in measure
from bounded inputs.
\end{proposition}

\begin{proof}
First let $f$ be bounded. Proposition~\ref{prop:mean} applies
to $f+r$ and $f-r$, and $R(f\pm r)=Rf$. Adding the resulting
inequalities and subtracting $2U_1(r,0)$ gives
\[
 \begin{split}
 &\int_X[U_1(f+r,Rf)+U_1(f-r,Rf)-2U_1(r,0)]\,d\mu\\
 &\hspace{12mm}\le
 U_1(\mathfrak m+r,0)+U_1(\mathfrak m-r,0)-2U_1(r,0).
 \end{split}
\]
Multiply by $c_pr^{p-2}$ and integrate over $r>0$.
The proof of Lemma~\ref{lem:superposition} gives
\[
 \begin{split}
 &c_p\int_0^\infty
 |U_1(x+r,y)+U_1(x-r,y)-2U_1(r,0)|r^{p-2}\,dr\\
 &\hspace{12mm}\le u_p(x,y)+a_p\le |x|^p+2a_p.
 \end{split}
\]
This bound is uniform in $y$ and integrable after substitution
$(x,y)=(f,Rf)$, so Fubini and \eqref{eq:superposition} prove
\eqref{eq:p-mean}. Since $B_p=U_p+|x|^p$, applying that estimate
to $f/\lambda$ proves \eqref{eq:B-finite} for bounded inputs.
In particular, \eqref{eq:B-bounds} gives
\[
 \mu\{|Rf|\ge\lambda\}\le K_p\lambda^{-p}\|f\|_p^p.
\]

Let bounded $f_n$ converge to $f$ in $L^p$. Applying this
distribution estimate to differences shows that $Rf_n$ is
Cauchy in measure. A subsequence chosen with summable
exceptional-set measures converges almost everywhere, and the
same estimate makes its limit independent of the approximants;
this is the extension argument of \cite[Lemma~3.11]{BHV}.
Choose the subsequence so that $f_n\to f$ almost everywhere as
well. Continuity and nonnegativity of $B_p$ permit Fatou in
\eqref{eq:B-finite}. The means converge to $\mathfrak m$, so the
mean term is retained. Finally, $U_p=B_p-|f|^p$ and
\eqref{eq:B-bounds} make \eqref{eq:p-mean} well defined and
recover it from the case $\lambda=1$.
\end{proof}

\begin{remark}\label{rem:compensation}
The compensation in \eqref{eq:superposition} must remain inside
the integral: its separate terms have divergent integrals in $r$.
The mean term in Proposition~\ref{prop:mean} supplies precisely
this compensation; the estimate $\int U_1(f,Rf)\le0$ alone does not.
At $p=2$,
\begin{equation}\label{eq:phi-p2}
 a_2=1,\qquad u_2(x,y)=x^2+1-y^2\quad(|y|<1),\qquad
 B_2(x,y)=\Phi_2(y)=\min\{y^2,1\}.
\end{equation}
Thus \eqref{eq:B-finite} becomes
\[
 \int_X\min\{|Rf|^2,\lambda^2\}\,d\mu
 \le\|f\|_2^2-|\mathfrak m|^2,
\]
consistent with $L^2$ contraction and $R\one=0$.
\end{remark}

\section{Passage to the classical transforms}
\label{sec:classical}

We approximate the transverse torus by finite conservative
generators and then let the whole torus grow. Conservativity is
needed for the prescribed-mean argument; compression of an
infinite transverse space may introduce killing.

\begin{proposition}\label{prop:classical-B}
Let $1<p\le2$, $n\ge1$, $1\le j\le n$, $\Omega=\T^n$ or $\R^n$, and
$f\in L^p(\Omega;\R)$. For every $\lambda>0$,
\begin{equation}\label{eq:classical-B}
 \lambda^p\int_\Omega
 B_p(f/\lambda,R_jf/\lambda)\,d\mu\le\|f\|_p^p.
\end{equation}
On the normalized torus, with
$\mathfrak m=\int_{\T^n}f\,d\mu$, the stronger estimate holds:
\begin{equation}\label{eq:torus-mean}
 \lambda^p\int_{\T^n}
 B_p(f/\lambda,R_jf/\lambda)\,d\mu
 \le\|f\|_p^p+\lambda^pU_p(\mathfrak m/\lambda,0).
\end{equation}
\end{proposition}

\begin{proof}
It suffices to take $j=1$ and $\lambda=1$.

\emph{The torus.}
For $n\ge2$, let $Y=\T^{n-1}$ with normalized measure and
$L=-\Delta_Y$. Choose refining finite measurable partitions of
$Y$, with positive weights, whose conditional expectations
$E_N$ converge strongly to the identity in $L^2(Y)$ and
$L^p(Y)$. These conditional expectations preserve the integral
and satisfy $E_N\one=\one$. For $r>0$, on $V_N=\ran E_N$ put
\[
 Q_{N,r}=E_Ne^{-rL}E_N\big|_{V_N},\qquad
 L_{N,r}=r^{-1}(I_{V_N}-Q_{N,r}).
\]
The matrix $Q_{N,r}$ is self-adjoint, positivity preserving,
preserves $\one$, and has spectrum in $[0,1]$. The expansion
\[
 e^{-tL_{N,r}}=e^{-t/r}\sum_{k=0}^\infty
                 \frac{(t/r)^k}{k!}Q_{N,r}^k
\]
shows that $e^{-tL_{N,r}}$ is a conservative Markov semigroup.
Proposition~\ref{prop:p-finite} therefore applies to $E_Nf$.

We use the spectral approximation from
\cite[Proposition~3.10]{BHV} with these conservative matrices.
Extending $L_{N,r}$ by zero on $V_N^\perp$ gives, at fixed $r$,
\[
 r^{-1}(E_N-E_Ne^{-rL}E_N)
 \longrightarrow L_r:=r^{-1}(I-e^{-rL})
\]
strongly, with uniformly bounded operator norms. At each
nonzero horizontal frequency, continuous functional calculus,
followed by dominated convergence in the Fourier sum, gives
\[
 R_{L_{N,r}}E_Nf\longrightarrow R_{L_r}f
 \quad\hbox{in }L^2(\T\times Y).
\]
The spectral theorem and $r^{-1}(1-e^{-rs})\to s$ similarly give
$R_{L_r}f\to R_Lf$ in $L^2$ as $r\downarrow0$.
For $f\in L^2(\T^n)\subset L^p(\T^n)$, take almost-everywhere
convergent subsequences at each step. Fatou's lemma and
$\|E_Nf\|_p\le\|f\|_p$ pass \eqref{eq:B-finite} through both
limits. Since $\int E_Nf\,d\mu=\mathfrak m$, this proves
\eqref{eq:torus-mean}, including its mean term, for $L^2$ inputs.
When $n=1$, take $E$ to be a single point and $L=0$.

\emph{Euclidean space.}
Let $f\in C_c^\infty(\R^n;\R)$. On the torus of side length
$\ell$, periodize $f$ to $f_\ell$ and denote its Riesz transform
by $R_{j,\ell}$. Simultaneous scaling of all coordinates and
multiplication by the torus volume give, for sufficiently large
$\ell$,
\begin{equation}\label{eq:large-torus}
 \int_{[-\ell/2,\ell/2)^n}B_p(f_\ell,R_{j,\ell}f_\ell)\,dx
 \le\|f\|_p^p.
\end{equation}
Here we have discarded the nonpositive mean term. On this
fundamental cube $f_\ell=f$ for large $\ell$. Moreover,
\[
 R_{j,\ell}f_\ell(x)
 =\ell^{-n}\sum_{k\in\mathbb Z^n}
 m_j(2\pi k/\ell)\widehat f(2\pi k/\ell)e^{2\pi ik\cdot x/\ell},
 \qquad m_j(\xi)=-i\xi_j/|\xi|,
\]
with $m_j(0)=0$. These Fourier series converge locally uniformly
to the Fourier integral for $R_jf$. Indeed, Riemann sums converge
uniformly on compact $x$-sets when the frequencies range over a
compact set away from zero. The sum and integral over
$|\xi|<\eta$ are bounded by $C(\eta+\ell^{-1})^n$, and the rapid
decay of $\widehat f$ controls the frequency tails uniformly for
$\ell\ge1$. Letting $\ell\to\infty$ and then $\eta\downarrow0$
proves the claim, as in the periodization argument of
\cite[Lemma~3.9]{BHV}, now in all coordinates.

Extend the nonnegative integrand of \eqref{eq:large-torus} by
zero outside its cube. Fatou's lemma proves
\eqref{eq:classical-B} for compactly supported smooth inputs.
On both spaces, the distribution estimate obtained from
\eqref{eq:B-bounds}, applied to differences of approximants,
and Fatou extend the conclusion to all real $L^p$ inputs,
as in Proposition~\ref{prop:p-finite}. On the torus their
means also converge, preserving \eqref{eq:torus-mean}.
The resulting transforms agree with the usual $L^p$ Riesz
transforms: they agree on the dense smooth class, whose
classical transforms converge in $L^p$ and hence in measure.
Finally, applying the estimates to $f/\lambda$ gives the
stated scaling.
\end{proof}

\begin{corollary}\label{cor:refinement}
Under the hypotheses of Proposition~\ref{prop:classical-B},
\begin{equation}\label{eq:phi-classical}
 \lambda^p\int_\Omega
 \Phi_p(R_jf/\lambda)\,d\mu\le\|f\|_p^p.
\end{equation}
Equivalently, the weak-type estimate has the nonnegative remainder
\begin{equation}\label{eq:phi-remainder}
 \begin{split}
 \lambda^p\mu\{|R_jf|\ge\lambda\}
 &+K_p\lambda^p\int_{\{|R_jf|<\lambda\}}
       \Phi_p(R_jf/\lambda)\,d\mu\\
 &\le K_p\|f\|_p^p.
 \end{split}
\end{equation}
On the normalized torus, in both inequalities the quantity
$\|f\|_p^p$ may be replaced by
$\|f\|_p^p+\lambda^pU_p(\mathfrak m/\lambda,0)$, where
$\mathfrak m=\int_{\T^n}f\,d\mu$.
\end{corollary}

\begin{proof}
The comparison $\Phi_p(y)\le B_p(x,y)$ and
Proposition~\ref{prop:classical-B} give \eqref{eq:phi-classical}
and its torus refinement. Since $\Phi_p=a_p$ on $\{|y|\ge1\}$
and $K_pa_p=1$, splitting the integral over the level set
and its complement gives \eqref{eq:phi-remainder}.
Dropping its nonnegative remainder proves the upper estimate
in Theorem~\ref{thm:main} for $1<p\le2$.
\end{proof}

\section{Sharpness}
\label{sec:sharpness}

We recall the one-dimensional extremal from
\cite[p.~910]{Janakiraman} in the circle form used in
\cite[Proposition~3.13]{BHV}, then transfer its lower bound to
each Euclidean dimension. The transfer uses the convergence of
the Riesz multiplier to the Hilbert-transform multiplier. Hence
we implement that mechanism by transverse dilation.

\begin{lemma}\label{lem:line-lower}
For $1<p\le2$, the optimal coefficient in the real-line Hilbert
transform distribution inequality is at least $K_p$.
\end{lemma}

\begin{proof}
On the normalized circle, consider the conjugate pair
\begin{equation}\label{eq:extremal}
 g_*(\theta)=\frac2\pi\log\left|\cot\frac\theta2\right|,
 \qquad Hg_*=\sgn(\sin\theta).
\end{equation}
The identity follows from the $L^2$ Fourier series
$g_*=(4/\pi)\sum_{k\ge0}\cos((2k+1)\theta)/(2k+1)$
and the sine series of $\sgn(\sin\theta)$.
The function $g_*$ has mean zero, belongs to every finite $L^p$, and
\begin{equation}\label{eq:extremal-norm}
 \|g_*\|_{L^p(\T)}^p
 =\frac1\pi\int_\R
       \frac{|(2/\pi)\log|t||^p}{1+t^2}\,dt=a_p.
\end{equation}
The substitution $t=\cot(\theta/2)$ proves the displayed equality.
Hence thresholds increasing to one give the lower coefficient
$a_p^{-1}$ on the circle. Real, mean-zero trigonometric
polynomial approximants may be used instead: convergence in
$L^2$ implies convergence in $L^p$ for $p\le2$ and convergence
of their Hilbert transforms in measure.

We record the passage to the line. Let $g$ be a real, mean-zero
trigonometric polynomial, extended periodically with period $2\pi$,
and set
\[
 F_N(x)=g(x)\one_{[-\pi N,\pi N]}(x),
 \qquad N\in\mathbb N.
\]
The periodic Hilbert transform equals the principal-value integral
of the periodic extension, with tails grouped in periods; this
follows termwise for $g$. Its zero mean supplies a bounded periodic
primitive $G$. Integration by parts on each omitted half-line gives
\begin{equation}\label{eq:tail}
 |H_{\R}F_N(x)-H_{\T}g(x)|
 \le\frac{C_g}{\pi N-|x|},\qquad |x|<\pi N.
\end{equation}
Both the endpoint term and the integral involving $G(t)/(x-t)^2$
are bounded by the reciprocal distance to the endpoint.

For $0<\alpha<1$, the error in \eqref{eq:tail} tends uniformly
to zero on $|x|\le(1-\alpha)\pi N$. If a coefficient $C$ is valid
for the line transform, apply it to $F_N$, using
\[
 \|F_N\|_{L^p(\R)}^p
 =2\pi N\|g\|_{L^p(\T)}^p.
\]
For $\lambda,\eta>0$, periodic averaging on the central interval
and \eqref{eq:tail} give, after division by $2\pi N$ and passage
to $N\to\infty$,
\[
 C\|g\|_p^p
 \ge\lambda^p(1-\alpha)
       m\{|H_\T g|>\lambda+\eta\}.
\]
Let $\alpha\downarrow0$ and $\eta\downarrow0$. Thus every
coefficient valid on the line is valid for trigonometric
polynomials on the circle. Approximation of \eqref{eq:extremal}
and thresholds increasing to one show $C\ge K_p$.
\end{proof}

\begin{proof}[Completion of the proof of Theorem~\ref{thm:main}]
The upper bound is Proposition~\ref{prop:classical-B}, together
with the quoted $p=1$ theorem. On $\T^n$, an input depending only
on the $j$th coordinate has $R_j$ equal to the one-dimensional
Hilbert transform. Equations \eqref{eq:extremal}--\eqref{eq:extremal-norm}
therefore give sharpness directly.

On $\R^n$ it suffices, by permutation of coordinates, to take
$j=1$. For $n\ge2$, choose a cube $Q\subset\R^{n-1}$ of volume
one and a real $g\in L^p(\R)\cap L^2(\R)$. Put
\[
 f_\varepsilon(x_1,x')=g(x_1)\one_Q(\varepsilon x'),
 \qquad
 V_\varepsilon(x_1,z')
 =(R_1f_\varepsilon)(x_1,z'/\varepsilon).
\]
After transverse rescaling the multiplier is
\[
 \frac{-i\xi_1}{(\xi_1^2+\varepsilon^2|\xi'|^2)^{1/2}}.
\]
It is bounded by one and converges almost everywhere to
$-i\sgn\xi_1$. Plancherel's theorem and dominated convergence give
\begin{equation}\label{eq:dilation}
 V_\varepsilon\longrightarrow
 (Hg)(x_1)\one_Q(z')\quad\hbox{in }L^2(\R^n).
\end{equation}
If $C$ is any valid coefficient for $R_1$ in dimension $n$,
rescaling its inequality gives
\[
 \lambda^p|\{|V_\varepsilon|>\lambda\}|
 \le C\|g\|_p^p,
\]
because both the original level-set measure and the input norm to
the $p$th power carry the factor $\varepsilon^{-(n-1)}$.
Pass along an almost-everywhere convergent subsequence in
\eqref{eq:dilation} and use Fatou for the open level sets.
It follows that
\[
 \lambda^p|\{|Hg|>\lambda\}|\le C\|g\|_p^p.
\]
The compactly supported inputs used in Lemma~\ref{lem:line-lower}
belong to $L^p\cap L^2$, so that lemma implies $C\ge K_p$.
For $n=1$ it applies directly. Sharpness at $p=1$ is already
part of \cite[Theorem~1.1]{BHV}.
\end{proof}

\begin{remark}[Equality in the torus integral estimates]\label{rem:equality}
For $1<p\le2$, take $f(x)=g_*(x_j)$ on $\T^n$. Then
$\mathfrak m=0$, $|R_jf|=1$ almost everywhere, and
$\|f\|_p^p=a_p$. Thus, at $\lambda=1$,
\[
 \int_{\T^n}B_p(f,R_jf)\,d\mu
 =\int_{\T^n}\Phi_p(R_jf)\,d\mu
 =a_p=\|f\|_p^p.
\]
This gives equality in \eqref{eq:classical-B},
\eqref{eq:torus-mean}, and both estimates of
Corollary~\ref{cor:refinement}. In the strict distribution
inequality \eqref{eq:main-weak}, the thresholds must still
approach one from below.
\end{remark}

\begin{remark}[The range $p>2$]\label{rem:p-greater-two}
Sharp weak-type inequalities for martingale transforms and
differentially subordinated martingales have been studied extensively,
starting with Burkholder's sharp weak-type $(1,1)$
inequality~\cite{Burkholder1979}. For real-valued martingales under
differential subordination, the optimal weak-type $(p,p)$ constants
were determined by Burkholder~\cite{Burkholder1984} for
$1\le p\le2$ and by Suh~\cite{Suh2005} for $2<p<\infty$.
For further sharp estimates under additional hypotheses, see
Os\k{e}kowski~\cite{Osekowski2009}, and for a systematic treatment
see \cite{Osekowski2012}.

To the best of our knowledge, the optimal
constant $C_p$ in the distribution inequality for the Hilbert transform
\begin{equation}\label{eq:hilbert-open}
 \lambda^p m\bigl(\{x\in\mathbb T:|Hf(x)|>\lambda\}\bigr)
 \le C_p\|f\|_{L^p(\mathbb T,m)}^p,
 \qquad \lambda>0,
\end{equation}
remains unknown for $2<p<\infty$; see
Ga{\l}\k{a}zka and
Os\k{e}kowski~\cite[Section~1]{GalazkaOsekowski2024}.
Here $H$ is the Hilbert transform on $\T$ and $m$ is
normalized Haar measure.

Ga{\l}\k{a}zka and Os\k{e}kowski~\cite[Theorem~1.2]{GalazkaOsekowski2024}
determine the optimal constant $A_{p,q}$ in the integral estimate
\begin{equation}\label{eq:hilbert-integral}
 \int_E |Hf|\,dm
 \le A_{p,q}\|f\|_{L^p(\mathbb T,m)}\,m(E)^{1-1/q},
 \qquad 1<q<p<\infty,
\end{equation}
for every real-valued $f\in L^p(\mathbb T,m)$ and every measurable
$E\subseteq\mathbb T$. The corresponding sharp estimate
with $q=p$, whose optimal constant we denote by $A_{p,p}$,
is recalled in their Theorem~1.1.
Applying this diagonal estimate to
$E=\{x\in\mathbb T:|Hf(x)|>\lambda\}$ gives
$C_p\le A_{p,p}^{\,p}$.
However, sharpness of the integral estimate does not imply
sharpness of the resulting distribution estimate, so these
results do not determine the optimal constant $C_p$ for $p>2$.

Similar sharp integral estimates for the individual Riesz
transforms on $\R^n$ were obtained by
Os\k{e}kowski~\cite[Theorem~1.1]{Osekowski2014}
for $1<p<\infty$. As in the case of the Hilbert transform,
sharpness in the integral formulation does not determine
the optimal constants in the usual distribution weak-type
$(p,p)$ inequalities for $2<p<\infty$.
\end{remark}

\end{document}